\documentclass[12pt]{amsart}
\usepackage{amssymb}
\usepackage{graphics,graphicx,color,amsmath,amssymb,latexsym}
\usepackage{epsbox}
\theoremstyle{plain}
\newtheorem{thm}{Theorem}[section]

\theoremstyle{definition}

\newtheorem{rmk}[thm]{Remark}

\def\pmc#1{\setbox0=\hbox{#1}
    \kern-.1em\copy0\kern-\wd0
    \kern.1em\copy0\kern-\wd0}

\def\De{\Delta}

\def\ka{\kappa}

\def\si{\sigma}

\def\vp{\varphi}
\def\om{\omega}
\def\op{\operatorname}
\def\ov{\overline}

\def\sm{\setminus}

\def\spmapright#1{\smash{%
\mathop{\hbox to 1.3cm {\rightarrowfill}}
\limits^{#1}}}
\def\sbmapright#1{\smash{%
\mathop{\hbox to 1.3cm {\rightarrowfill}}
\limits_{#1}}}

\begin{document}
\title
{From uncountable abelian groups to uncountable nonabelian groups}
\author{Katsuya Eda}
\address{Department of Mathematics, 
Waseda University, Tokyo 169-8555, JAPAN}
\email{eda@waseda.jp}
\begin{abstract}
The present note surveys my research related to generalizing notions of abelian group theory to non-commutative group case and applying them particularly to the investigation of fundamental groups. 
\end{abstract}
\subjclass[2010]{20K20, 20K25, 55N10, 55Q52}
\keywords{one-dimensional, Peano continuum, fundamental group, singular homology, Hawaiian earring, 
algebraic compactness, cotorsion, shape group}

\maketitle
\section{Introduction}
This paper is an account on my studies of topics in mathematics and, although they are rooted in abelian group theory, they mostly only indirectly are related to abelian groups themselves. 
The emphasis is to show connections between my study of abelian groups to that of fundamental groups, which are non-abelian.  
To state theorems exactly we need to use technical terms from algebraic topology, for which we refer the reader to \cite{Spanier:algtop}. But, I take care so that the reader can understand the outline without understanding precise definitions, which is the main content of this paper. 
\section{The Specker theorem}\label{sec:specker}
My joining into abelian group people started from my attendance of the Honolulu conference in 1982-1983. Before then, R\"{u}diger G\"{o}bel, who passed away in 2014, contacted me as one of the organizers of the conference. Laszlo, who by then already was a central person in abelian group theory, was present so that I met him there. After then I worked on abelian groups for several years. The reason why I started my study of abelian groups is my interest to the Specker theorem about $\mathbb{Z}^\om$ \cite{Specker:ZN}, which is a unique theorem about infinitely generated discrete groups supporting a duality theorem, i.e. 
\[
\op{Hom}(\op{Hom}(\oplus _\om \mathbb{Z}, \mathbb{Z}), \mathbb{Z}) \cong \oplus _\om \mathbb{Z}. 
\]
I felt that there should exist good mathematics around this theorem. Since the \v{C}ech homology group 
$\check{H}_1(\mathbb{H})$ of the Hawaiian earring $\mathbb{H}$ is isomorphic to $\mathbb{Z}^\om$, I aimed at analyzing the algebraic structure of the Hawaiian earring (see the figure).

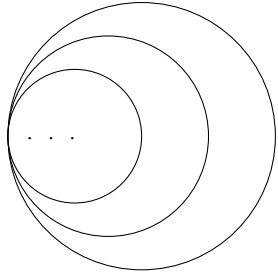
\begin{center}
  \psbox[scale=0.7]{gask3.eps}
\end{center}

I tried to find applications of the Specker theorem and also the Chase lemma to algebraic topology. 
Since the singular homology group of the Hawaiian earring is not isomorphic to $\mathbb{Z}^\om$, contrasting to the \v{C}ech homology group, I introduced a canonical factor of the singular homology \cite{ES:singular} and investigated such groups for spaces of certain types. This will be explained in Section 3. But, I felt it is not sufficient as what I felt from the Specker theorem and tried to investigate the fundamental group of the Hawaiian earring and to prove a non-commutative version of the Specker theorem. 
After having obtained a proof  around 1985 I only found out that G. Higman's work \cite{Higman:unrestrict} from 1952 already contained a demonstration of this fact: Namely, there he first proves that every 
 homomorphism from the unrestricted free product, i.e. the canonical inverse limit of free groups of finite rank, factors through a free group of finite rank with the projection, which 
can be seen a non-commutative version of the Specker theorem. 
Then, Higman mentions the validity of this result also for a certain subgroup $P$ of the inverse limit. 
Three years later, H. B. Griffiths \cite{Griffiths:product} proved $P$ to be isomorphic to $\pi _1(\mathbb{H})$. 
I introduced a new notion free $\si$-product $\pmc{$\times$}\, \, \;_{i\in I}^\si G_i$ 
of groups $G_i$ to investigate fundamental groups of spaces like 
the Hawaiian earring \cite{E:free}. A free $\si$-product $\pmc{$\times$}\, \, \;_{i\in I}^\si G_i$ is a subgroup 
of the unrestricted free product in \cite{Higman:unrestrict} consisting of elements expressed by words 
defined on countable linearly ordered sets, while an element of a usual free product is expressed by words defined on finite linearly ordered sets. 
There, I proved a non-commutative version of the Chase lemma, which I've mentioned above, i.e.  
\begin{thm}\label{thm:chase1}\cite[Theorem 2.1]{E:free} Let 
\[
h: \pmc{$\times$}\, \, \;_{i\in I}^\si G_i \to \ast _{j\in J}H_j
\]
be a homomorphism for groups $G_i$ and $H_j$. Then, there exist finite subsets $F$ of $I$ and 
$G$ of $J$ respectively such that 
\[
h(\pmc{$\times$}\, \, \;_{i\in I\sm F}^\si G_i) \le \ast _{j\in G}H_j.
\]
\end{thm}
The original Chase lemma is about homomorphisms from direct products to direct sums. In the category of abelian groups there exist injective objects and consequently the statement is more complicated, but in the non-abelian case it becomes simpler. 
Actually the conclusion of Theorem~\ref{thm:chase1} can be considerably sharpened as we will see in Theorem~\ref{thm:chase2}. 

Applying Higman's theorem I proved that every endomorphism on $\pi _1(\mathbb{H})$ is conjugate to an endomorphism induced from a continuous self-map on $\mathbb{H}$. Though it has become a seminal result, I felt that it was very far from results about the Sierpinski carpet and Menger sponge at that time, since the Hawaiian earring has only one wild point, while all the points are wild in the others (see the sequel to Theorem~\ref{thm2} for the word "wild"). 

At that time, conjugators about the fundamental groups were troublesome for me. But, a few years later I found that conjugators become the key for finding points 
of the above spaces in their fundamental groups, fundamental groupoids more exactly. 
Actually we have: 
\begin{thm}\cite[Theorem 1.1]{E:spatial} 
Let $X$ be a one-dimensional metric space and $h:\pi _1(\mathbb{H},o)\to \pi _1(X,x_0)$ a homomorphism. 
Then, there exists a continuous map $f:\mathbb{H}\to X$ and 
a point $x\in X$ and a path $p$ from $x_0$ to $x$ such that $f(o)=x$ and 
$h = \vp _p\circ f_*$, where $\vp _p: \pi _1(X,x)\to \pi _1(X,x_0)$ is the point change isomorphism and 
$f_*$ is the homomorphism induced by $f$. If the image of $h$ is uncountable, the point $x$ is unique and 
$p$ is unique up to homotopy relative to end points. 
\end{thm}
A conjugator in the fundamental group, which was troublesome, corresponds to the homotopy type of a loop $p$ in the statement. Points can be restored from fundamental groups as the maximal compactible families of subgroups which are the homomorphic images of $\pi _1(\mathbb{H})$ \cite{ConnerEda:information}. 

Based on this we have,  
\begin{thm}\label{thm2}\cite[Theorem 1.3]{E:spatial} 
Let $X$ and $Y$ be one-dimensional, locally path-connected, path-connected metric spaces which are not semi-locally simply connected at any point. If the fundamental groups are isomorphic, then $X$ and $Y$ are homeomorphic. 
\end{thm}
A one-dimensional space space is called {\it semi-locally simply connected}, if any point has a neighborhood without a circle. We call a space {\it wild}, if the space contains a point at which the space 
is not semi-locally simply-connected. 
Theorem~\ref{thm2} implies that the fundamental groups of the Sierpinski carpet and the Menger sponge are not isomorphic to each other, since the two spaces are not homeomorphic. This result was quite unexpectable: 
The homotopy equivalence of spaces implies the isomorphic-ness of fundamental groups but in general the homotopy equivalence is much weaker than the homeomorphism type. Poincar\'{e} introduced the notion of fundamental groups, as a much rough equivalence in comparison with a homeomorphism type. Therefore, though it is a very restricted case, this was unexpected. Although I several times had already checked my proof of Theorem~\ref{thm2}, I still distrusted it and consequently even tried to manufacture a counter example at least a few times, and I've heard that several topologists did not believe Theorem~\ref{thm2}. 

Then, I worked on this line, i.e. to investigate relationship between properties of groups and those of spaces. This duality between spaces and groups through the fundamental groups is extended to at most countable direct products and other constructions \cite{ConnerEda:information} of spaces. 
But it was difficult to publish such papers, since no one except me was working in this area. Since my retirement year was 2017, I needed to publish my papers as a researcher and worked on other subjects. Among them I proved: 
\begin{thm}\label{thm:3}\cite[Theorem 1.1]{E:onedim} 
For one dimensional Peano continua $X$ and $Y$, $X$ and $Y$ are homotopy equivalent, if and only if  
$\pi _1(X)$ and $\pi _1(Y)$ are isomorphic. 
\end{thm}
Although this result looks like a standard statement in algebraic topology, but the proof of this theorem depends on Theorem~\ref{thm2}, which is extraordinary. Consequently this theorem is actually an extraordinary theorem. 

Before then, as I mentioned, Theorem~\ref{thm:chase1} was strengthened as follows. 
\begin{thm}\label{thm:chase2}\cite[Theorem 1.3]{E:atom} 
Under the same assumption of Theorem~\ref{thm:chase1},  there exist a finite subset $F$ of $I$ and an element $j\in J$ such that 
$h(\pmc{$\times$}\, \, \;_{i\in I\sm F}^\si G_i)$ is contained in a subgroup conjugate to $H_j$. 
\end{thm}
As a variant of this theorem, we have 
\begin{thm}\label{thm:chase3}\cite[Theorem 1.4]{E:atom} 
Let $X$ be a path-connected, locally path-connected, first countable space which is not semi-locally simply-connected at any point. If $h: \pi _1(X)\to \ast _{j\in J}H_j$ is an injective homomorphism, then 
the image of $h$ is contained in a subgroup conjugate to some $H_j$. 
\end{thm}

These are results where I translated known facts from abelian group theory to not necessarily abelian goups, and, in particular, to fundamental groups. 

\section{Algebraically compact groups}\label{sec:algcompact}
As is well-known, the singular homology group of a path-connected space is the abelianizations of the fundamental group for a path-connected space. 
In algebraic topology it is well-known that all groups appear as fundamental groups and, consequently, all abelian groups appear as homology groups. However, the corresponding spaces constructed by using group theoretic data are artificial ones. 
On the other hand, the fundamental groups of spaces which have local complexities, e.g. fractals, are 
out of studies for a long time. Also, divisible, or algebraically compact groups did not occur as homology groups or their subgroups of spaces, which are familiar as topological spaces, i.e. not artificial or formal 
objects. 

Back to around 1990, I found out that the singular homology groups of certain spaces 
are {\it complete modulo the Ulm subgroup} \cite{E:union} (a notion which had been  
introduced by Dugas-Goebel \cite{DugasGoebel:algcompact}). In particular, 
the divisible group $\mathbb{Q}$ occurs as a subgroup.  
Since torsion-free groups which are complete modulo the Ulm subgroup are algebraically compact, 
reduced algebraically compact groups possibly occur as subgroups of the singular homology 
group of the Hawaiian earring at that time \cite{E:free}, see Figure. 

Actually the first singular homology group of the Hawaiian earring, i.e. $H_1(\mathbb{H})$, turned out to be isomorphic to 
\[
\mathbb{Z}^\om \oplus \oplus _{\mathfrak{c}}\mathbb{Q} \oplus \Pi _{p:\rm{prime}}A_p, 
\]
where $\mathfrak{c}$ is the cardinality of the continuum and $A_p$ is the $p$-completion of the free abelian group of rank $\mathfrak{c}$ \cite{EK:hawaii}.  
After proving Theorem~\ref{thm:3}, I investigated the singular homology groups of one dimensional Peano continua. My assessment was that he singular homology group of the Hawaiian earring can hardly isomorphic to any of the Menger sponge or the Sierpinski carpet - because of the so much simpler topological structure of the Hawaiian earring. 
In spite of my many trials to show this, I had failed to do it. After changing my mind I proved, 
\begin{thm}\cite{E:singular}
The singular homology group $H_1(X)$ of a one-dimension Peano continuum $X$ is isomorphic to a free abelian group of finite rank or the singular homology group $H_1(\mathbb{H})$ of the Hawaiian earring. 
\end{thm}
That is, if a Peano continuum $X$ contains a single wild point, then $H_1(X)$ is already isomorphic to $H_1(\mathbb{H})$, no matter how many additional wild points are there in $X$.  
\begin{rmk}
In \cite{E:free} I used the notion ''complete modulo the Ulm subgroup'' due to 
Dugas-Goebel \cite{DugasGoebel:algcompact}. Very recently Herfort-Hojka \cite{HerfortHojka:cotorsion} have characterized the notion ''cotorsion'' using equation systems. According to it results in \cite{E:free} 
can be improved to being cotorsion. This characterization is a commutative version of the equation system due to G. Higman \cite{Higman:unrestrict} for non-commutative groups. 
\end{rmk}

\section{The Reid class and \v{C}ech systems}  
The Reid class consists of the integer group $\mathbb{Z}$ and the groups obtained by iterating use of 
forming direct sums and direct products \cite{Reid:almostfreegroup}. Back to 1983, answering a question in \cite{DugasZimmermann:sumproduct, Zimmermann-Huisgen:Fuchs76}, whether $\op{Hom}(A,\mathbb{Z})$ belongs to the Reid class for every abelian group $A$, I proved that $\op{Hom}(C(\mathbb{Q},\mathbb{Z}), \mathbb{Z})$ does not belong to the Reid class, where $C(\mathbb{Q},\mathbb{Z})$ is a group of continuous functions. 
Before then I was interested in algebraic topology and so I tried to apply the Reid class to algebraic topology. There are two typical and distinct ways of attaching infinitely many circles with a common point. The one is so-called a bouquet, where a basic open neighborhood consists of open neighborhoods of all circles. A basic open neighborhood of the common point in the other way consists of almost all copies of the circle and open neighborhoods in the remaining finite circles. When the number of copies are countable, 
the latter space is homeomorphic to the Hawaiian earring, while the former one is called a countable bouquet.  Passing to the factor group of singular homology as described in \cite{ES:singular} one obtains 
respectively free abelian groups $\oplus _\ka \mathbb{Z}$ and direct products $\mathbb{Z}^\kappa$. Let us recall this factorization \cite{ES:singular}: 
The singular chain group $S_n(X)$ is the free abelian group generating by the set $C(\De _n,X)$, where $\De _n$ is the $n$-simplex. Since $C(\De _n,X)$ has the compact open topology, $S_n(X)$ can be regarded as the free abelian topological group over $C(\De _n,X)$. Then, the boundary operator $\partial _n:S_{n+1}(X)\to S_n(X)$ becomes continuous. Therefore, 
$\op{Ker}(\partial _{n-1})$ is closed, but $\op{Im}(\partial _n)$ may not be closed. We take its closure 
and consider $\op{Ker}(\partial _n)/\ov{\op{Im}}(\partial _{n+1})= \ov{H}_n(X)$, which is the factor of singular homology. 
By taking the closure, information due to the wildness of a space disappears. For instance, 
$H^T_1(\mathbb{H})$ is isomorphic to $\mathbb{Z}^\om$. 
On the other hand, if $X$ is locally good, e.g. locally contractible, $\op{Im}(\partial _{n+1})$ is closed and we have $H^T_n(X)=H_n(X)$. Therefore, for spaces usually appearing in algebraic topology it gives us the same as singular homology. 

When attaching infinitely many component spaces with one common point 
we have two typical types of one point unions.  A neighborhood of the common point is a union of neighborhoods of the point in component space in one type and is a union of neighborhoods of the point in finitely many component 
spaces and the whole spaces for remaining component spaces. 
These construction can be done alternately and iterated. In such constructions with the same common point 
the complexity of the topology around the common point should increase. The complexity can be expressed by the hierarchy theorem of the Reid class via the factor of singular homology \cite{ES:singular} and \cite{E:Z-kernel, EklofMekler:almostfree} (I talked about this in the Oberwolfach conference in 1989). 

Many years later I, with J. Nakamura, tried to classify the inverse limits of sequences of finitely generated free groups.
Such inverse limits in the abelian case become finitely generated free abelian groups or $\mathbb{Z}^\om$, which is a consequence of the fact that any subgroup of a finitely generated free abelian group is also  finitely generated. Remark that a subgroup of a finitely generated free group may not be finitely generated. 
 The inverse limits of sequences of finitely generated free groups are precisely 
the first \v{C}ech homotopy groups of one-dimensional connected compact metric spaces. When such spaces are locally connected, they are Peano continua. Then, we have an inverse system of surjective homomorphisms and the inverse limits are isomorphic to finitely generated free groups $F_n$ or the inverse limit according to the canonical projections of  finitely generated free groups. 
That is, $G_1 =\displaystyle{ \varprojlim _{n\to\infty} \ast _{i=0}^n\mathbb{Z}_i}$, where the bonding map from $\ast _{i=0}^n\mathbb{Z}_i$ to $\ast _{i=0}^{n-1}\mathbb{Z}_i$ is the projection. 
In general, we have three other groups. Let $F_\om$ be the countable free group. 
Let $G_2 =\displaystyle{ \varprojlim _{n\to\infty} F_\om \ast\ast _{i=0}^n\mathbb{Z}_i}$, where the bonding map from $F_\om \ast\ast _{i=0}^n\mathbb{Z}_i$ to $F_\om \ast\ast _{i=0}^{n-1}\mathbb{Z}_i$ is the projection. Let $F_{\om ,n}$ be copies of $F_\om$ and 
$G_3 = 
\displaystyle{ \varprojlim _{n\to\infty} \ast _{i=0}^n F_{\om ,i}}$ where the bonding map from $\ast _{i=0}^n F_{\om ,i}$ to $\ast _{i=0}^{n-1}F_{\om ,i}$ is the projection, 
 is the inverse limit according to the projections of finite free products of copies of $F_\om$.
Now, $F_\om$, $G_2$ and $G_3$ are the three groups. 
 As you can see from these results, the classification of the inverse limit of at most countable free groups is the same as the corresponding one to finitely generated free groups \cite{EN:inverse}. Therefore, the classification of shape groups of one dimensional connected, compact metric spaces , i.e. their first \v{C}ech homotopy groups, is the same as that of the corresponding groups of connected separable metric spaces. 

We explain how we proved that these uncountable groups $G_1, G_2$ and $G_3$ are not isomorphic to each others. 
Let $[G,G]$ be the commutator subgroup of a group $G$ and $\op{Ab}(G) = G/[G,G]$. Let 
$R_\mathbb{Z}(A) = \bigcap \{ \op{Ker}(h)\, | \, h\in \op{Hom}(A,\mathbb{Z})\}$. As is well-known, $A/R_\mathbb{Z}(A)$ is isomorphic to a subgroup of $\mathbb{Z}^\ka$ for some $\ka$.  
We consider a functor $\op{F}(G) =\op{Ab}(G)/R_\mathbb{Z}(\op{Ab}(G))$. 
Then, we have $\op{F}(G_1) = \mathbb{Z}^\om , \op{F}(G_2) = \oplus _\om \mathbb{Z}\oplus \mathbb{Z}^\om$ and $\op{F}(G_3) = (\oplus _\om \mathbb{Z})^\om$. Then, by the hierarchy theorem of the Reid class we conclude these groups are not isomorphic and consequently so are $G_1, G_2$ and $G_3$.

Together with the results in Sections~\ref{sec:specker} and \ref{sec:algcompact}, I feel the following. Though abelian groups are related to many areas of mathematics through homology and cohomology, apart from the finitely generated case the relationships are formal ones. Therefore, what I've explained in the preceding are new aspects in the relationships between infinitely generated abelian groups, non-abelian groups and topological spaces. 
In particular, divisible subgroups or algebraically compact (cotorsion) subgroups of singular homology groups imply the existence of wild points in the spaces. 

I am still working on the topic discussed here and so please do not consider this note as my final report on these issues.  
I hope I can continue my research until the age of Laszlo, i.e. twenty years more. 
When I wrote to him "I understand that one can continue doing mathematics over 90,"
he shot back a proposal to argue that doing mathematics actively would get me over 90.  
\medskip

\begin{center}
{\bf Acknowledgement}
\end{center}
The author thanks the referee for improving wording to express contents exactly and reading this manuscript thoroughly. It seems that the referee also read some old references. 


\end{document}